\documentclass[11pt]{article}
\usepackage{amsfonts}
\usepackage{latexsym,amsmath}
\usepackage{amssymb,array}
\makeatletter \oddsidemargin 0in \evensidemargin 0in \textwidth
16cm \RequirePackage[dvips]{graphicx} \textheight 20cm
\newtheorem{t1}{Theorem}[section]
\newtheorem{p1}{Proposition}[section]
\newtheorem{l1}{Lemma}[section]
\newtheorem{c1}{Corollary}[section]
\newtheorem{d1}{Definition}[section]
\newtheorem{r1}{Remark}[section]

\newtheorem{ex}{Example}[section]

\begin{document}
\title{On Cumulative Residual (Past) Extropy of Extreme Order Statistics}
\author{Chanchal Kundu\footnote{Corresponding
author e-mail: chanchal$_{-}$kundu@yahoo.com, ckundu@rgipt.ac.in.}\\
Department of Mathematics\\
Rajiv Gandhi Institute of Petroleum Technology\\
Jais 229 304, U.P., India}
\date{April, 2020}
\maketitle
\begin{abstract}
In the recent information-theoretic literature, the concept of extropy has been studied for order statistics. In the present communication we consider a cumulative analogue of extropy in the same vein of cumulative residual (past) entropy and study it in context with extreme order statistics. A dynamic version of cumulative residual (past) extropy for smallest (largest) order statistic is also studied here. It is shown that the proposed measures (and their dynamic versions) of extreme order statistics determine the distribution uniquely. Some characterizations of the generalized Pareto and power distributions, which are commonly used in reliability modeling, are given.
\end{abstract}
\section{Introduction}
The idea of information-theoretic entropy was introduced by Shannon (1948) which plays an important role in diverse areas such as financial analysis, data compression, molecular biology, hydrology, meteorology, computer science and information theory. In spite of enormous success of Shannon's entropy, this measure has some drawbacks and may not be appropriate in every situation. To get rid of these drawbacks an alternate measure of uncertainty called cumulative residual (past) entropy (CREn) (CPEn) has been introduced in the literature that extends Shannon's entropy.\\
\hspace*{.2in} Let $X$ be an absolutely continuous nonnegative random variable with probability density function ({\it pdf}) $f$, cumulative distribution function ({\it cdf}) $F(x)$ and survival function ({\it sf}) $\overline F(x)$. Rao et al. (2004) defined CREn as
\begin{eqnarray}\label{eq1}
\varepsilon(X)=-\int_{0}^{\infty}\overline F(x)\ln{\overline
F(x)}dx.
\end{eqnarray}
The idea was to replace the {\it pdf} by {\it sf} in Shannon's differential entropy
\begin{eqnarray*}
H(X)=-\int_{0}^\infty f(x)\ln f(x) dx.
\end{eqnarray*}
For more properties and applications of CREn one may refer to Rao (2005). Later, Di Crescenzo and Longobardi (2009) introduce a dual measure based on the
{\it cdf} $F(x)$, called the cumulative
past entropy (CPEn) which is analogous to CREn, as follows:
\begin{eqnarray}\label{eq2}
\overline\varepsilon(X)=-\int_{0}^{\infty}F(x)\ln{F(x)}dx.
\end{eqnarray}
Several aspects of dynamic CREn (CPEn) has been studied by Asadi and Zohrevand (2007), Navarro et al. (2010), Khorashadizadeh et al. (2013) and Di Crescenzo and Longobardi (2013), among others.\\
\hspace*{.2in} The notion of entropy is recently entwined with a complementary dual measure, designated as extropy, by Lad et al. (2015). With this concept, they provide a completion to theories of information based on entropy, resolving a longstanding question in its axiomatization as proposed by Shannon and pursued by Jaynes. The extropy of the random variable $X$ is defined as
\begin{eqnarray}\label{eq3} J(X)=-\frac{1}{2}\int_0^\infty f^2(x)dx.
\end{eqnarray}
As is entropy, the extropy is interpreted as a measure of the amount of uncertainty represented by $X$. However, the two measures are distinct and fundamentally intertwined with each other. The discovery of extropy was stimulated by a problem that arises in the application of the theory of proper scoring rules for alternative forecast distributions. One statistical application of the concept was provided to the scoring of sequential forecast distributions. As argued by Lad et al. (2015), in any commercial and scientific arena in which entropic computations have become standard, its complementary dual would also be well worth investigating. Recently, Qiu and Jia (2018a) consider a dynamic version of extropy and study its properties. Qiu and Jia (2018b) proposed two estimators of extropy and develop a goodness-of-fit test for standard uniform distribution. For some bounds on extropy with variational distance one may refer to Yang et al. (2018). Further, Qiu et al. (2018) explored some extropy properties of mixed systems.\\
\hspace*{.2in} Recently, Jahanshahi et al. (2019) introduce a cumulative analogue of extropy, called cumulative residual extropy (CREx) which is analogous to CREn, as follows:
\begin{eqnarray}\label{eq4}\varepsilon J(X)=-\frac{1}{2}\int_0^\infty \overline F^2(x)dx.
\end{eqnarray}
 This measure is defined for lifetime distributions, taking into account the {\it sf} $\overline F(x)$ instead of the {\it pdf} $f$ as has been done in CREn. It gives a chance for the {\it sf} to be involved directly in the measurement of extropy. Recall that the equilibrium random variable of $X$, denoted by $X_e$, has the density function $f_e(x)=\overline F(x)/E(x),~x\in(0,1)$ where $E(X)$ is finite. It deserves interest in reliability and queueing theory. Clearly, $J(X_e)=\varepsilon J(X)/(E(X))^2$ giving that CREx is, apart from a constant term, a measure of extropy of $X_e$ in the unity measure of $E(X)$. For some properties, estimations and applications of CREx one may refer to Jahanshahi et al. (2019). They have successfully applied CREx in risk measure and to measure the independence between two random variables. \\
\hspace*{.2in} Along a similar line, the cumulative past extropy (CPEx) can be defined as
\begin{eqnarray}\label{eq5}\overline\varepsilon J(X)=-\frac{1}{2}\int_0^\infty F^2(x)dx.
\end{eqnarray}
The basic idea is to replace the density function by distribution function in extropy (\ref{eq3}) as the distribution function exists even in cases where density does not. Moreover, in practice what is of interest and/or measurable is the distribution function. For example, if the random variable is the life span of a machine, then the event of interest is not whether the life span equals $t$, but rather whether the life span exceeds $t$. With this, the investigation of CPEx would be well worth investigating in contrast with CREx.\\
\hspace*{.2in} Let $X_1, X_2,\ldots,X_n$ be $n$ independent and identically distributed (iid) random variables having a univariate continuous {\it cdf} $F$. The order statistics of the sample are defined by the arrangement of
$X_1, X_2,\ldots,X_n$ from the smallest to the largest, denoted as $X_{1:n}\leqslant X_{2:n}\leqslant\ldots\leqslant X_{n:n}$. These statistics have been used in a wide range of problems, including robust statistical estimation, detection of outliers, characterization of probability distributions and goodness-of-fit tests, analysis of censored samples, reliability analysis, quality control and strength of materials, see, for details, Arnold et al. (1992), David and Nagaraja (2003), and references therein. The minimum and maximum are examples of extreme order statistics and are defined by
$$X_{1:n}={\rm min}\{X_1, X_2,\ldots,X_n\}~~{\rm and}~~X_{n:n}={\rm max}\{X_1, X_2,\ldots,X_n\}.$$
The extrema $X_{1:n}$ and $X_{n:n}$ are of special interest in many practical problems of distributional analysis. The extremes arise in the statistical study of floods and droughts, as well as in problems of breaking strength and fatigue failure. Recently, Qiu (2017) studied extropy of order statistics and record values exploring some properties and characterization results. Most recently, Qiu and Jia (2018a) provided some results on residual extropy of order statistics. For the study on CREn/CPEn of order statistics one may refer to Baratpour (2010), Thapliyal et al. (2013) and Park and Kim (2014), to mention a few. But to the best of our knowledge, no attention has been paid to the study of cumulative extropy of order statistics. Motivated by this, in this article, we consider CREx (CPEx) for smallest (largest) order statistics.\\
\hspace*{.2in} The rest of the paper is organized as follows. In Section 2 we consider several properties of CREx of $X_{1:n}$ and obtain some bounds. It is shown that the CREx of $X_{1:n}$ uniquely determines the parent distribution. An analogous discussion is made on dynamic CREx of $X_{1:n}$ and some orderings for an $(n-k+1)-$out$-$of$-n$ system are obtained. Generalized Pareto distribution is also characterized. Section 3 is devoted to the study of similar results for CPEx of $X_{n:n}$. Some characterizations of power distribution are also provided here. Finally, in Section 4, a conclusion is made on the present study.
\section{Results on (dynamic) CREx of smallest order statistic}
In this section first we review some properties of (dynamic) CREx for smallest order statistic. We also define a stochastic ordering based on DCREx. Later, we provide characterization of some well-known distributions including uniqueness of DCREx for smallest order statistic.
\subsection{Some properties}
Let $X_{1:n}$ be the smallest order statistic in a random sample of size $n$ from an absolutely continuous nonnegative random variable $X$. Then the distribution function of $X_{1:n}$ is given by $F_{X_{1:n}}(x)=1-\overline F^n(x).$
Thus, \begin{eqnarray}\label{eq2.1}
\varepsilon J(X_{1:n})=\frac{-1}{2}\int_0^\infty\overline F^{2n}(x)dx.
\end{eqnarray}
By change of variable $u=\overline F(x)$, we have
\begin{eqnarray}\label{eq2.2}
\varepsilon J(X_{1:n})=\frac{-1}{2}\int_0^1 \frac{u^{2n}}{f\left(F^{-1}(1-u)\right)}du.
\end{eqnarray}
Expression of the measure (\ref{eq2.1}) for some specific distributions are given below.
\begin{ex}\label{ex2.1}
\begin{itemize}
\item [(i)] If $X$ is uniformly distributed on $[a,b],~a<b$, then $\varepsilon J(X_{1:n})=\frac{-(b-a)}{2(2n+1)}$.
\item [(ii)] Let $X$ follow finite range distribution $\overline F(x)=(1-ax)^b;~a,b>0$ and $x\in(0,1/a)$, then $\varepsilon J(X_{1:n})=\frac{-1}{2a(1+2nb)}$.
\item [(iii)] For Weibull distribution with $\overline F(x)=e^{-\lambda x^\theta},~\lambda,\theta>0$, $\varepsilon J(X_{1:n})=\frac{-\Gamma\left(\frac{1}{\theta}\right)}{2\theta[2n\lambda ]^{1/\theta}}$.
\item [(iv)] If $X$ follows Folded Cramer distribution $\overline F(x)=\frac{1}{1+\theta x},~x,\theta>0$ then $\varepsilon J(X_{1:n})=\frac{-1}{2(2n-1)\theta}$.
\item [(v)] For Pareto distribution $\overline F(x)=\left(\frac{\lambda}{x+\lambda}\right)^\theta,~x,\lambda>0$ and $\theta>1$, $\varepsilon J(X_{1:n})=\frac{-\lambda}{2(2n\theta-1)}$.
\end{itemize}
\end{ex}
\hspace*{.2in} First, let us look at the monotonic behavior of (\ref{eq2.1}). The following lemma will be useful.
\begin{l1}\label{lm1} $h(\alpha)=u^\alpha$ is a decreasing function of $\alpha$ for $u\in(0,1)$.
\end{l1}
\begin{p1}\label{p2.1} For $n\geqslant 1$, $\varepsilon J(X_{1:n})$ is increasing in $n$.
\end{p1}
Proof: On using Lemma \ref{lm1}, it can easily be checked from (\ref{eq2.2}) that $$\varepsilon J(X_{1:n})\geqslant\varepsilon J(X_{1:n-1}).$$
Hence the result follows. $\hfill\square$\\

\hspace*{.2in} Now, we examine some lower bounds of $\varepsilon J(X_{1:n})$ (cf. Jahanshahi et al., 2019). The proof follows from Lemma \ref{lm1} and the above proposition.
\begin{p1}\label{p2.2} Let $X_{1:n}$ be the smallest order statistic in an iid random sample $X_1,X_2,\ldots,X_n$ from a nonnegative absolutely continuous random variable $X$ with finite mean $\mu$. Then
\begin{itemize}
  \item [i.] $\varepsilon J(X_{1:n})\geqslant\frac{-\mu}{2}$;
  \item [ii.] $\varepsilon J(X_{1:n})\geqslant\varepsilon J(X)$.
\end{itemize}
The second result elucidates that the uncertainty of $X$ is less than that of $X_{1:n}$ and it is in contrast with the similar investigation based on CREn.
\end{p1}
\hspace*{.2in} The following lemma popularly known as M$\ddot{u}$ntz-Sz$\acute{a}$sz theorem, which is often invoked in moment based characterizations (see Kamps, 1998) will be used to prove the upcoming theorems.
\begin{l1}\label{lm2} For any strictly increasing sequence of positive integers $\{n_j,~j\geqslant1\}$, the sequence of polynomials $\{x^{n_j}\}$ is complete on $L(0,1)$ if and only if $\sum_{j=1}^\infty n_j^{-1}=+\infty.$
\end{l1}
\hspace*{.2in} In the sequel we assume that $\{n_j,~j\geqslant1\}$ is a strictly increasing sequence of positive integers. The following theorem shows that the parent distribution can be characterized by CREx of $X_{1:n}$.
\begin{t1}\label{th2.1} Let $X$ and $Y$ be nonnegative absolutely continuous random variables with distribution functions $F(x)$ and $G(x)$, respectively. Then $F$ and $G$ belong to the same family of distributions, but for a change in location, if and only if
\begin{eqnarray}\label{eq2.3}\varepsilon J(X_{1:n})=\varepsilon J(Y_{1:n}),\end{eqnarray}
for all $n=n_j,~j\geqslant1$ such that $\sum_{j=1}^\infty n_j^{-1}=+\infty.$
\end{t1}
Proof: The necessity is trivial. For the sufficient part let (\ref{eq2.3}) holds, then we have
$$\int_0^1u^{2n}\left[\frac{1}{f\left(F^{-1}(1-u)\right)}-\frac{1}{g\left(G^{-1}(1-u)\right)}\right]du=0.$$
The rest of the proof is similar to Theorem 2.2 of Baratpour (2010). $\hfill\square$
\begin{c1}Taking $n=1$ in the above theorem it can easily be shown that CREx, $\varepsilon J(X)$ uniquely determines the distribution function.
\end{c1}
\hspace*{.2in} In the following theorem we obtain a characterization result for the scale family of distributions.
\begin{t1}\label{th2.2}Let $X$ and $Y$ be two random variables as described in Theorem \ref{th2.1}. Also let $X$ and $Y$ have common support $[0,\infty)$. Then $F$ and $G$ belong to the same family of distributions, but for a change of scale, if and only if
\begin{eqnarray}\label{eq2.4}\varepsilon J(X_{1:n})/E(X_{1:n})=\varepsilon J(Y_{1:n})/E(Y_{1:n}),\end{eqnarray}
for all $n=n_j,~j\geqslant1$ such that $\sum_{j=1}^\infty n_j^{-1}=+\infty.$
\end{t1}
Proof: The necessary part is trivial. To prove the sufficient part suppose (\ref{eq2.3}) holds, then on noting that
$$E(X_{1:n})=\int_0^{\infty}\overline F^n(x)dx=\int_0^1u^n/f\left(F^{-1}(1-u)\right)du$$
we have
$$\frac{\int_0^1\frac{u^{2n}}{f\left(F^{-1}(1-u)\right)}du}{\int_0^1\frac{u^{2n}}{g\left(G^{-1}(1-u)\right)}du}
=\frac{\int_0^1\frac{u^n}{f\left(F^{-1}(1-u)\right)}du}{\int_0^1\frac{u^{n}}{g\left(G^{-1}(1-u)\right)}du}=\frac{1}{c},~{\rm say}.$$
Or equivalently,
\begin{eqnarray}\label{eq2.3+1}\int_0^1u^{2n}\left[\frac{1}{f\left(F^{-1}(1-u)\right)}-\frac{c^{-1}}{g\left(G^{-1}(1-u)\right)}\right]du=0.
\end{eqnarray}
If (\ref{eq2.3+1}) holds for $n=n_j,~j\geqslant1$, such that $\sum_{j=1}^\infty n_j^{-1}=+\infty,$ then from Lemma \ref{lm2} we can conclude that $f\left(F^{-1}(w)\right)=cg\left(G^{-1}(w)\right)$ for all $0<w<1$. Which on using the fact $\frac{d}{dw}F^{-1}(w)=\frac{1}{f\left(F^{-1}(w)\right)}$, reduces to $F^{-1}(w)=(1/c)G^{-1}(w)+d$ where $c(>0)$ and $d$ are arbitrary constants. Since $X$ and $Y$ have common support $[0,\infty)$, therefore $d$ must be equal to zero. This means $F$ and $G$ belong to the same family of distributions, but for a change in scale. $\hfill\square$\\

\begin{figure}
\centering
\includegraphics[height=4cm,keepaspectratio]{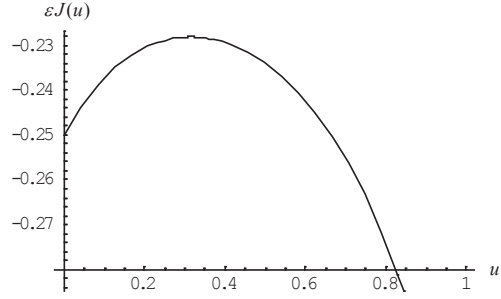}
\caption{Plot of $\varepsilon J(u)$ for $u\in(0,1)$.}\label{fig2.1}
\end{figure}

\hspace*{.2in} Consider the following example to verify the above result.
\begin{ex} Let $X$ and $Y$ follow exp($\lambda_1$) and exp($\lambda_2$), respectively with common support $[0,\infty)$. Here $X$ and $Y$ belong to the same family but having different scale. It can easily be seen that $\varepsilon J(X_{1:n})/E(X_{1:n})=-1/4=\varepsilon J(Y_{1:n})/E(Y_{1:n})$. Converse part follows along the lines of the proof of the theorem.
\end{ex}
\hspace*{.2in} The study of duration is a subject of interest in many areas such as reliability, survival analysis, actuary, economics, business etc. Capturing effects of the age $t$ of an individual or an item under study on the information about the remaining lifetime is important in many applications and so the dynamic (time dependent) information measures have been considered in the literature. Let $X$ be the lifetime of a component or system under the condition that the system has survived to age $t$. In such a case, the distribution of interest, for computing uncertainty and information, is the residual lifetime $X_t=(X-t|X>t)$. The dynamic CREx (DCREx) for the residual lifetime distribution $X_t$ is defined as
$$\varepsilon J(X;t)=\frac{-1}{2}\int_t^\infty\left(\frac{\overline F(x)}{\overline F(t)}\right)^{2}dx.$$
Analogously, the DCREx  of $X_{1:n}$ is given
\begin{eqnarray}\label{eq2.5}
\varepsilon J(X_{1:n};t)=\frac{-1}{2}\int_t^\infty\left(\frac{\overline F(x)}{\overline F(t)}\right)^{2n}dx.
\end{eqnarray}
In the sequel we study the monotonic behavior of $\varepsilon J(X_{1:n};t)$. Differentiating (\ref{eq2.5}) with respect to $t$, we get
\begin{eqnarray}\label{eq2.6}
\frac{d}{dt}\varepsilon J(X_{1:n};t)=2n\lambda_F(t)\varepsilon J(X_{1:n};t)+1/2,
\end{eqnarray}
where $\lambda_F(t)=f(t)/\overline F(t)$ is the hazard rate of $X$. Therefore, $\varepsilon J(X_{1:n};t)$ is increasing (decreasing) in $t$, if and only if
$$\varepsilon J(X_{1:n};t)\geqslant(\leqslant)\frac{-1}{4n\lambda_F(t)}.$$
It is to be mentioned here that exponential distribution is the borderline distribution as seen in most of the cases. For finite range distribution DCREx is increasing in $t$ while for Pareto distribution it is decreasing. To see that not all distributions are monotone in terms of DCREx, we consider a random variable having survival function $\overline F(x)=1-(1-e^{-x})(1-e^{-2x})$, $x>0$. Then, Figure \ref{fig2.1} shows that DCREx is not monotone. Note that the substitution $t=-\log u$ has been used while plotting curve so that $\varepsilon J(X;t)=\varepsilon J(u)$, say.\\
%--------------------------------------
%-----------------------------------------

\hspace*{.2in} In the following we give some lower bounds of $\varepsilon J(X_{1:n};t)$ in a similar line of Proposition \ref{p2.2}.
\begin{p1}\label{p2.3} Let $X_{1:n}$ be the smallest order statistic in an iid random sample $X_1,X_2,\ldots,X_n$ from a nonnegative absolutely continuous random variable $X$ with mean residual life (mrl) $\delta_F(t)=E(X_t)$. Then
\begin{itemize}
  \item [(i)] $\varepsilon J(X_{1:n};t)\geqslant\frac{-\delta_F(t)}{2};$
  \item [(ii)] $\varepsilon J(X_{1:n};t)$ is increasing in $n$ and further $\varepsilon J(X_{1:n};t)\geqslant\varepsilon J(X;t)$.
\end{itemize}
\end{p1}
\hspace*{.2in} The following example illustrate the above theorem.
\begin{ex} Let $X$ follow finite range distribution as given in Example \ref{ex2.1}, then
$$\varepsilon J(X_{1:n};t)=-\left(\frac{1+b}{1+2nb}\right)\frac{\delta_F(t)}{2}\geqslant\frac{-\delta_F(t)}{2}.$$
Thus, (i) holds. On the other hand,
$$\varepsilon J(X_{1:n};t)-\varepsilon J(X;t)=\frac{2b(n-1)(1-at)}{2a(1+2nb)(1+2b)}\geqslant0,~\forall n\geqslant1,$$
which confirms (ii). $\hfill\square$
\end{ex}

\hspace*{.2in} Now we discuss a ordering based on DCREx. For more details on stochastic orders one may refer to Shaked and Shanthikumar (2007).
\begin{d1} The random variable $X$ is said to be less than or equal to $Y$ in the dynamic CREx, denoted by $X\leqslant_{DCREx}Y$, if $\varepsilon J(X;t)\geqslant\varepsilon J(Y;t).$
\end{d1}
\begin{ex} Suppose that $X_i$ follows Pareto distribution $\overline F_i(x)=\left(\frac{\lambda}{x+\lambda}\right)^{\theta_i},~x,\lambda>0$ and $\theta_i>1,~i=1,2$. Now, in view of $\varepsilon J(X_i;t)=-(\lambda+t)/(4\theta_i-2)$, one can easily see that if $\theta_1>\theta_2$ then $X_1\leqslant_{DCREx}X_2$.
\end{ex}
\begin{r1} It is not very difficult to check that the DCREx order is a partial order i.e., reflexive, anti-symmetric and transitive properties are satisfied for this order.
\end{r1}
\hspace*{.2in} The following theorem shows that the DCREx order between two nonnegative random variables is closed under increasing linear transformation. The proof being straightforward is omitted. It is worthwhile to remark that the monotonicity property of DCREx is also preserved under such transformation.
\begin{t1} For two nonnegative random variables $X_1$ and $X_2$, let $Y_i=aX_i+b,~i=1,2$ with $a>0$ and $b\geqslant0$. Then $Y_1\leqslant_{DCREx}Y_2$ provided $X_1\leqslant_{DCREx}X_2$.
\end{t1}
\hspace*{.2in} Now we have the more general result. For brevity, the proof is omitted.
\begin{t1} For two nonnegative random variables $X_1$ and $X_2$, let us define $Y_i=a_iX_i+b_i,~i=1,2$ where $a_1\geqslant a_2>0$ and $b_1\geqslant b_2>0$. Then $Y_1\leqslant_{DCREx}Y_2$ if $X_1\leqslant_{DCREx}X_2$ and either $\varepsilon J(X_1;t)$ or $\varepsilon J(X_2;t)$ is decreasing in $t$.
\end{t1}

\hspace*{.2in} The next theorem is related to the DCREx and hazard rate (hr) orderings. The proof follows from (\ref{eq2.5}).
\begin{t1}\label{th2.4} Let $X$ and $Y$ be two absolutely continuous nonnegative random variables with hazard functions $\lambda_F,~\lambda_G$, respectively. If $X\leqslant_{hr}Y$, i.e., $\lambda_F(t)\geqslant\lambda_G(t)$ for all $t\geqslant0$, then $X_{1:n}\leqslant_{DCREx}Y_{1:n}$.
\end{t1}
\hspace*{.2in} An application of the above theorem is given in the following example.
\begin{ex} If $X$ follows Weibull distribution with shape parameter $\theta$, then for $\theta_1>\theta_2$, $X_{\theta_1}\leqslant_{hr}X_{\theta_2}$. Therefore, $X_{1:n;\theta_1}\leqslant_{DCREx}X_{1:n;\theta_2}$. Note that for any random variable $X$, $\lambda_{F_{1:n}}(t)\geqslant\lambda_F(t)$ and thus $\varepsilon J(X_{1:n};t)\geqslant\varepsilon J(X;t)$, as shown in Proposition \ref{p2.3}.
\end{ex}
\hspace*{.2in} Let $X_{k:n}$ be the $k^{th}$ order statistic in a set of iid random variables from $F$. It represents the lifetime of an $(n-k+1)-$out$-$of$-n$ system which is a common structure of redundancy and widely used in reliability theory and survival analysis. It is to be noted that $X_{1:n}$ represents the lifetime of a series system, whereas $X_{n:n}$ that of a parallel system. The following result is due to Nagaraja (1990) which will be used to prove the upcoming theorem.
\begin{t1}\label{th2.5} Let $X_{k_1:n_1}$ and $X_{k_2:n_2}$ be two order statistics from two sets of iid random variables having common distribution $F$ with sizes $n_1$ and $n_2$, respectively. Then $\lambda_{F_{k_2:n_2}}(t)=\theta(t)\lambda_{F_{k_1:n_1}}(t)$ such that $\theta(t)$ is increasing in $t$ with $\theta(t)\in(0,1)$ for
\begin{itemize}
  \item $n_1=n_2=n,~k_1=k,~k_2=k+1;$
  \item $n_1=n,~n_2=n-1,~k_1=k_2=k;$
  \item $n_1=n,~n_2=n+1,~k_1=k,~k_2=k+1.$
\end{itemize}
\end{t1}
\hspace*{.2in} Now on using Theorems \ref{th2.4} and \ref{th2.5} we have the following result for comparison of generalized order statistics.
\begin{t1}\label{th2.6} Let $X_{k:n}$ be the $k^{th}$ order statistic in a set of iid random variables from $F$. Then
\begin{itemize}
  \item $\varepsilon J(X_{k:n};t)\geqslant\varepsilon J(X_{k+1:n};t),$
  \item $\varepsilon J(X_{k:n};t)\geqslant\varepsilon J(X_{k:n-1};t),$
  \item $\varepsilon J(X_{k:n};t)\geqslant\varepsilon J(X_{k+1:n+1};t).$
\end{itemize}
\end{t1}
\subsection{Characterizations}
In the literature, the problem of characterizing probability distributions has been investigated by many researchers. The standard practice in modeling statistical data is either to derive the appropriate model based on the physical properties of the system or to choose a flexible family of distributions and then find a member of the family that is appropriate to the data. In both the situations it would be helpful if we find characterization theorems that explain the distribution. In fact, characterization approach is very appealing to both theoreticians and applied workers. The general characterization problem is to obtain whether $\varepsilon J(X_{1:n};t)$ uniquely determines the underlying distribution function. We consider the following characterization result in analogy with Theorem \ref{th2.1}.
\begin{t1} Let $X$ and $Y$ be nonnegative absolutely continuous random variables with distribution functions $F(x)$ and $G(x)$, respectively. Then $F$ and $G$ belong to the same family of distributions, but for a change in location, if and only if
\begin{eqnarray*}
\varepsilon J(X_{1:n};t)=\varepsilon J(Y_{1:n};t),\end{eqnarray*}
for all $n=n_j,~j\geqslant1$ such that $\sum_{j=1}^\infty n_j^{-1}=+\infty.$
\end{t1}
\begin{r1}For $n=1$, the above theorem gives that DCREx, $\varepsilon J(X;t)$ characterizes the
distribution function of $X$ uniquely, which can be proved otherwise also on similar lines.
\end{r1}
\hspace*{.2in} In reliability theory, in studies of the lifetime of a component or a system, a flexible model which has been
widely used in the literature is that of a generalized Pareto distribution (GPD) with {\it sf}
\begin{eqnarray}\label{eq2.7}
\overline F(x)=\left(\frac{\theta}{\lambda x+\theta}\right)^{\frac{1}{\lambda}+1},~x>0,~\theta>0,~\lambda>-1.
\end{eqnarray}
It plays an important role in extreme value theory and other branches of statistics. For the use of GPD in minimizing power consumption one may refer to Shang et al. (2003). The GPD, as a family of distributions, includes the exponential distribution when $\lambda\rightarrow0$, the Pareto Type-II distribution or Lomax distribution for $\lambda>0$, which is used in the studies of income, sizes of human settlements, reliability modeling and so on. The GPD becomes power distribution for $-1<\lambda<0$. In the following theorems, we obtain some results characterizing the GPD based on CREx of $X_{1:n}$.
\begin{t1}\label{th2.8}Let $X$ be a nonnegative absolutely continuous random variable with mrl $\delta_F(t)$. Then
\begin{eqnarray}\label{eq2.8}\varepsilon J(X_{1:n};t)=-c\delta_F(t),
\end{eqnarray}
where $c$ is a constant, characterizes GPD with survival function (\ref{eq2.7}). In particular, $X$ has
\begin{itemize}
  \item[(i)] exponential distribution iff $c=1/2n$,
  \item[(ii)] Pareto-II distribution iff $c<1/2n$,
  \item[(iii)] power distribution iff $1/2n<c<1$.
\end{itemize}
\end{t1}
Proof: For GPD (\ref{eq2.7}), (\ref{eq2.8}) is straightforward from (\ref{eq2.5}) with $c=1/[2(2n(1+\lambda)-\lambda)]$. Conversely, let (\ref{eq2.8}) holds. Then
$$\int_t^\infty \overline F^{2n}(x)dx=c\overline F^{2n-1}(t)\int_t^\infty \overline F(x)dx.$$
Differentiating both sides with respect to $t$, and using the relation $\lambda_F(t)\delta_F(t)=1+\delta'_F(t)$, we get
$$\delta'_F(t)=\frac{1-c}{(2n-1)c}-1=c_1-1,~{\rm say}.$$
This shows that mrl function of $X$ is linear. Thus, the desired result follows along with the distributions $(i)-(iii)$ according as $c_1=1,~c_1>1$ and $0<c_1<1$, respectively (cf. Hall and Wellner, 1981). $\hfill\square$\\

\hspace*{.2in} The following theorem gives another characterization of GPD.
\begin{t1}\label{th2.9} For an absolutely continuous nonnegative random variable $X$
\begin{eqnarray}\label{eq2.9}
\varepsilon J'(X_{1:n};t)=c,
\end{eqnarray}
where $c$ is a constant, if and only if $X$ follows GPD (\ref{eq2.7}).
\end{t1}
Proof: The if part follows from (\ref{eq2.5}) on noting that $\varepsilon J(X_{1:n};t)=-(\theta+\lambda t)/[2(2n(1+\lambda)-\lambda)]$. To prove the converse, let us assume that (\ref{eq2.9}) holds. Then, on using (\ref{eq2.6}), we obtain
$$-\frac{\lambda_F'(t)}{\lambda_F^2(t)}=\frac{4nc}{2c-1},$$
which on integration yields
$$\lambda_F(t)=\frac{1}{c_1t+c_2},$$
where $c_1=4nc/(2c-1)$ and $c_2^{-1}=\lambda_F(0)$. This is the hazard rate of GPD and the result follows by noting that the distribution function is determined uniquely by its hazard rate.
\section{(Dynamic) CPEx of largest order statistic}
In this section we introduce the concept of cumulative past extropy (CPEx) and an analogous discussion to Section 2 is made for CPEx of last order statistic $X_{n:n}$. To begin with, we investigate some basic properties of CPEx of a random variable $X$ having bounded support $[0,b]$ with $b$ finite. For the random variable $X$, we define the CPEx as
\begin{eqnarray}\label{eq3.1}
\overline\varepsilon J(X)=\frac{-1}{2}\int_0^bF^2(x)dx.
\end{eqnarray}
Outwardly, CPEx is a dual concept of CREx which relates extropy on the future lifetime of a system. It is suitable to measure information when extropy is related to the past. Let us analyze the effect of linear transformation on CPEx.
\begin{p1}\label{p3.1} Let $Y$ be a nonnegative random variable where $Y=cX+d$ with $c>0$ and $d\geqslant0$. Then $\overline\varepsilon J(Y)=c\overline\varepsilon J(X)$, which shows that CPEx is a shift-independent measure. Needless to say that this is reminiscent of the same property of CREn, CPEn and CREx.
\end{p1}
\hspace*{.2in} In the following, a relation is obtained between CPEx and CPEn which follows from inequality $-\log x\geqslant (1-x)$, for $x>0$.
\begin{p1} For the random variable $X$, $\overline\varepsilon J(X)\leqslant\frac{1}{2}\left[\overline\varepsilon(X)-(b-E(X))\right].$
\end{p1}
\hspace*{.2in} The following theorem shows that CPEx preserves the following mathematical property of CREn and CPEn. Note that the corresponding inequality for CREx as given in Jahanshahi et al. (2019) is quite different.
\begin{t1}\label{th3.1} Let $X$ and $Y$ be two nonnegative and independent random variables with supports $(0,b_X)$ and $(0,b_Y)$, respectively. Then $\overline\varepsilon J(X+Y)\geqslant\max\{\overline\varepsilon J(X),\overline\varepsilon J(Y)\}$.
\end{t1}
Proof: Since $X$ and $Y$ are independent, the {\it cdf} of $X+Y$ is given by $$F_{X+Y}(x)=\int_0^{b_Y}F_X(x-c)dF_Y(c),$$
where $F_Z(\cdot)$ is the {\it cdf} of $Z$. Using Jensen's inequality
$$F^2_{X+Y}(x)\leqslant\int_0^{b_Y}F_X^2(x-c)dF_Y(c).$$
Integrating both the sides with respect to $x$, we get
\begin{eqnarray*}
\overline\varepsilon J(X+Y)&\geqslant&\frac{-1}{2}\int_0^{b_X}\left(\int_0^{b_Y}F_X^2(x-c)dF_Y(c)\right)dx\\
&=&\frac{-1}{2}\int_0^{b_Y}dF_Y(c)\int_c^{b_X}F_X^2(x-c)dx\\
&=&\overline\varepsilon J(X),
\end{eqnarray*}
where the first equality follows on using the fact that $F_X(x-c)=0$ for $x\leqslant c$. Similarly, one can see that $\overline\varepsilon J(X+Y)\geqslant\overline\varepsilon J(Y)$. Hence the result follows. $\hfill\square$\\

\hspace*{.2in} In the following, we will find two more inequalities for CPEx.
\begin{t1} Let $X$ and $Y$ be two iid random variables with support $[0,b]$. Then
\begin{itemize}
  \item $E(|X-Y|)\geqslant4\overline\varepsilon J(X)$;
  \item $\overline\varepsilon J(X)\geqslant\left(\frac{E(X)-b}{2}\right)$.
\end{itemize}
\end{t1}
Proof: Since $X$ and $Y$ are iid random variables, then
$$2F(x)-2F^2(x)=P[\max(X,Y)>x]-P[\min(X,Y)>x].$$
Integrating both the sides we get
\begin{eqnarray}\label{eq3.2}
2\int_0^bF(x)\left(1-F(x)\right)dx=E[\max(X,Y)-\min(X,Y)]=E\left[|X-Y|\right].
\end{eqnarray}
The rest of the proof follows from (\ref{eq3.2}). $\hfill\square$\\

\hspace*{.2in} The following theorem gives a relationship between the conditional and the unconditional CPEx. It states that conditioning has a decreasing effect on CPEx.
\begin{t1}If $X$ and $Y$ are nonnegative random variables with supports $(0,b_X)$ and $(0,b_Y)$, respectively then $\overline\varepsilon J(X)\geqslant E_Y\left[\overline\varepsilon J(X|Y)\right]$.
\end{t1}
Proof: By using Jensen's inequality, it is not very difficult to see that
\begin{eqnarray*}E_Y\left[\overline\varepsilon J(X|Y)\right]&=&-\frac{1}{2}\int_0^{b_Y}\left(\int_0^{b_X}F^2_{X|Y}(x|y)dx\right)f_Y(y)dy\\
&\leqslant&-\frac{1}{2}\int_0^{b_X}\left(\int_0^{b_Y}F_{X|Y}(x|y)f_Y(y)dy\right)^2dx.
\end{eqnarray*}
Hence the result is obtained on noting that $\int_0^{b_Y}F_{X|Y}(x|y)f_Y(y)dy=F_X(x)$. $\hfill\square$\\

\hspace*{.2in} Let $X_{n:n}$ be the largest order statistic in a random sample of size $n$ from an absolutely continuous nonnegative random variable $X$. Then the distribution function of $X_{n:n}$ is given by $F_{X_{n:n}}(x)=F^n(x).$
Thus, \begin{eqnarray}\label{eq3.3}
\overline\varepsilon J(X_{n:n})=\frac{-1}{2}\int_0^b F^{2n}(x)dx.
\end{eqnarray}
For $n=1$ in (\ref{eq3.3}), we get (\ref{eq5}). In the following we study the monotonic behavior and bounds of $\overline\varepsilon J(X_{n:n})$. The proof is analogous to Proposition \ref{p2.2} and hence omitted.
\begin{p1} Let $X_{n:n}$ be the largest order statistic in an iid random sample $X_1,X_2,\ldots,X_n$ from an absolutely continuous random variable $X$ with finite support $[0,b]$ and mean $\mu$. Then
\begin{itemize}
  \item $\overline\varepsilon J(X_{n:n})\geqslant\frac{-1}{2}(b-\mu)$;
  \item $\overline\varepsilon J(X_{n:n})$ is increasing in $n$ and further $\overline\varepsilon J(X_{n:n})\geqslant\overline\varepsilon J(X)$.
\end{itemize}
\end{p1}

\hspace*{.2in} In the sequel we have some characterization results for location and scale family of distributions.
\begin{t1}\label{th3.4} Let $X$ and $Y$ be nonnegative absolutely continuous random variables with distribution functions $F(x)$ and $G(x)$, respectively. Then $F$ and $G$ belong to the same family of distributions, but for a change in location and scale, if and only if
\begin{eqnarray}\label{eq3.4}\overline\varepsilon J(X_{n:n})/\overline\varepsilon J(X)=\overline\varepsilon J(Y_{n:n})/\overline\varepsilon J(Y),
\end{eqnarray}
for all $n=n_j,~j\geqslant1$ such that $\sum_{j=1}^\infty n_j^{-1}=+\infty.$
\end{t1}
Proof: The necessary part is trivial. To prove the sufficient part suppose (\ref{eq3.4}) holds. Then, on using $v=F(x)$ in (\ref{eq3.3}), we have
$$\frac{\int_0^1\frac{v^{2n}}{f\left(F^{-1}(v)\right)}dv}{\int_0^1\frac{v^{2n}}{g\left(G^{-1}(v)\right)}dv}
=\frac{\int_0^1\frac{v^2}{f\left(F^{-1}(v)\right)}dv}{\int_0^1\frac{v^{2}}{g\left(G^{-1}(v)\right)}dv}=c,~{\rm say}.$$
Or equivalently,
\begin{eqnarray}\label{eq3.5}\int_0^1v^{2n}\left[\frac{1}{f\left(F^{-1}(v)\right)}-\frac{c}{g\left(G^{-1}(v)\right)}\right]dv=0.
\end{eqnarray}
If (\ref{eq3.5}) holds for $n=n_j,~j\geqslant1$, such that $\sum_{j=1}^\infty n_j^{-1}=+\infty,$ then from Lemma \ref{lm2} we can conclude that $1/f\left(F^{-1}(v)\right)=c/g\left(G^{-1}(v)\right)$ for all $0<v<1$, which on using the fact $\frac{d}{dv}F^{-1}(v)=\frac{1}{f\left(F^{-1}(v)\right)}$, reduces to $F^{-1}(v)=cG^{-1}(v)+d.$ This means $F$ and $G$ belong to the same family of distributions, but for a change in location and scale. $\hfill\square$\\

\hspace*{.2in} We consider the following example in support of the above theorem.
\begin{ex} Let $X$ and $Y$ belong to the same family of uniform distributions having different location and scale with cdfs $F(x)=(x-a_1)/(b_1-a_1),~0<a_1<x<b_1$ and $G(x)=(x-a_2)/(b_2-a_2),~0<a_2<x<b_2$, respectively. Then, straightforward calculation yields $\overline\varepsilon J(X_{n:n})/\overline\varepsilon J(X)=3/(2n+1)=\overline\varepsilon J(Y_{n:n})/\overline\varepsilon J(Y)$. We omit the converse part for brevity.
\end{ex}
\hspace*{.2in} In the following theorem we show that the parent distribution can be characterized by CPEx of $X_{n:n}$. The proof being similar to Theorem \ref{th2.1} is omitted.
\begin{t1}\label{th3.5} Under the assumptions of Theorem \ref{th3.4}, $F$ and $G$ belong to the same family of distributions, but for a change in location, if and only if
\begin{eqnarray*}\overline\varepsilon J(X_{n:n})=\overline\varepsilon J(Y_{n:n}),\end{eqnarray*}
for all $n=n_j,~j\geqslant1$ such that $\sum_{j=1}^\infty n_j^{-1}=+\infty.$
\end{t1}
\hspace*{.2in} Measure of uncertainty in past lifetime distribution plays an important role in the context of information theory, forensic sciences, and other related fields. For instance, if a system that begins to work at time $t=0$ is observed only at deterministic inspection times, and has been found inactive at time $t(>0)$, then $X_{[t]}=(t-X|X\leqslant t)$ represents the inactivity time of the system, i.e., the time elapsing between the inspection time $t$ and the failure time $X$. The dynamic CPEx (DCPEx) for the past lifetime distribution (or inactivity time) $X_{[t]}$ is defined as $$\overline\varepsilon J(X;t)=\frac{-1}{2}\int_0^t \left(\frac{F(x)}{F(t)}\right)^{2}dx.$$
It is not a shift-independent measure as for a random variable $Y$ as in Proposition \ref{p3.1}, $\overline\varepsilon J(Y;t)=c\overline\varepsilon J\left(X;\frac{t-d}{c}\right),~t>d$. If $X$ is symmetric with respect to $b/2$, i.e., $F(x)=\overline F(b-x)$ for $0\leqslant x\leqslant b$, then $\overline\varepsilon J(X;t)=\varepsilon J(X;b-t)$. Note that for lifetime distributions DCPEx is always negative and, decreasing in $t$. This decreasing property of $\overline\varepsilon J(X;t)$, on linear transformation of the random variable $X$, is preserved. The following example shows that DCPEx could be non-monotone as well.
\begin{ex}\label{ex3.2} Let $X$ be a nonnegative random variable with cdf
\begin{equation*}F(x)=\left\{\begin{array}{ll}\exp\{-1/2-1/x\},
\quad 0<x\leqslant 1\\
\exp\{-2+x^2/2\},
\quad 1<x\leqslant 2\\
1, \quad x\geqslant 2.
\end{array}\right.
\end{equation*}
%--------------------------------------
\begin{figure}
\centering
\includegraphics[height=4cm,keepaspectratio]{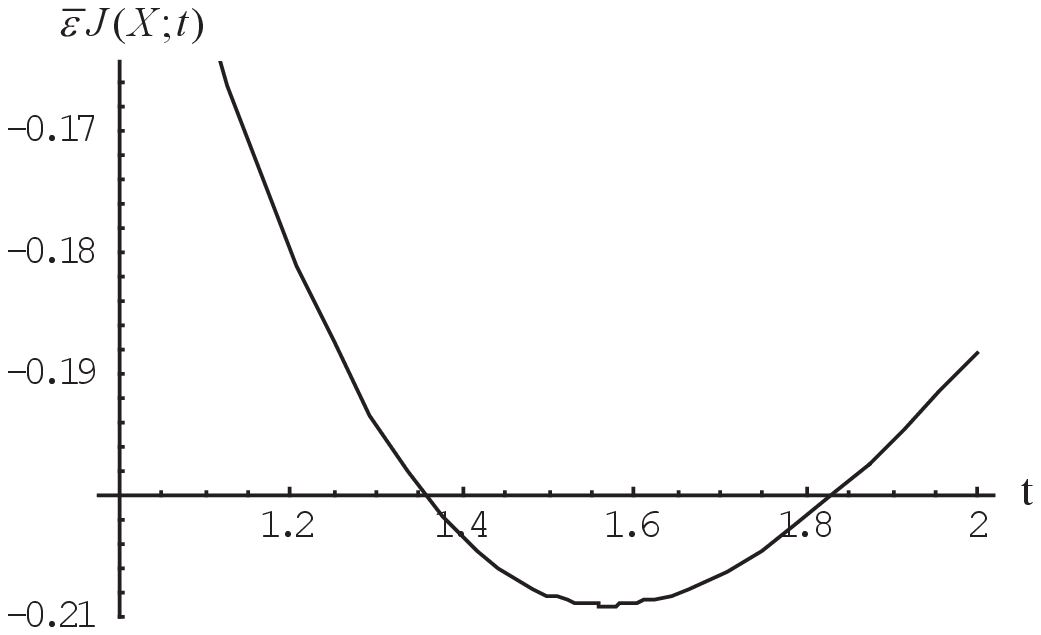}
\caption{Plot of $\overline\varepsilon J(X;t)$ for $t\in(1,2)$ (Example \ref{ex3.2}).}\label{fig3.1}
\end{figure}
%-----------------------------------------
Then $\overline\varepsilon J(X;t)$ is not monotone as shown in Figure \ref{fig3.1}.
\end{ex}
\hspace*{.2in} The DCPEx of $X_{n:n}$ is given by
\begin{eqnarray}\label{eq3.6}
\overline\varepsilon J(X_{n:n};t)=\frac{-1}{2}\int_0^t\left(\frac{F(x)}{F(t)}\right)^{2n}dx.
\end{eqnarray}
In analogy with Proposition \ref{p2.3}, below we discuss monotonic behavior and bounds of $\overline\varepsilon J(X_{n:n};t)$.
\begin{p1}\label{p3.4}Let $X_{n:n}$ be the largest order statistic in an iid random sample $X_1,X_2,\ldots,X_n$ from a nonnegative absolutely continuous random variable $X$ with expected inactivity time $m_F(t)=E(X_{[t]})$. Then
\begin{itemize}
    \item $\overline\varepsilon J(X_{n:n};t)$ is increasing in $n$ and decreasing in $t$;
  \item $\overline\varepsilon J(X_{n:n};t)\geqslant\frac{-1}{2} m_F(t)$;
  \item $\overline\varepsilon J(X_{n:n};t)\geqslant\overline\varepsilon J(X;t)$.
\end{itemize}
\end{p1}
\hspace*{.2in} Now we discuss a ordering based on DCPEx.
\begin{d1} For two random variables $X$ and $Y$, $X$ is said to be greater than $Y$ in DCPEx order (written as $X\geqslant_{DCPEx}Y$) if $\overline\varepsilon J(X;t)\geqslant\overline\varepsilon J(Y;t)$.
\end{d1}
\begin{ex}Let $X$ and $Y$ follow power distribution with cdfs $F(x)=x^{\alpha_1}$ and $G(x)=x^{\alpha_2}$ where $0<x<1$ and $\alpha_1,\alpha_2>0$. If $\alpha_1>\alpha_2$, then $X\geqslant_{DCPEx}Y$.
\end{ex}
\begin{r1} It can easily be verified that DCPEx order is also a partial order.
\end{r1}
\hspace*{.2in} The following theorem is immediate from (\ref{eq3.6}).
\begin{t1} Let $X$ and $Y$ be two absolutely continuous nonnegative random variables with reversed hazard rate functions $\phi_F,~\phi_G$ respectively. Then $\overline\varepsilon J(X_{n:n};t)\geqslant\overline\varepsilon J(Y_{n:n};t)$ if $X\geqslant_{rh}Y$ i.e., $\phi_F(t)\geqslant\phi_G(t)$ for all $t\geqslant0$.
\end{t1}
An immediate application of the theorem is that for a parallel system $\phi_{F_{n:n}}(t)=n\phi_F(t)\geqslant\phi_F(t)$. Therefore, $\overline\varepsilon J(X;t)\leqslant\overline\varepsilon J(X_{n:n};t)$ as also shown in Proposition \ref{p3.4}.\\
\hspace*{.2in} On using Theorem 2.2 and 2.4 of Kundu et al. (2009) we have the following result for an $(n-k+1)-$out$-$of$-n$ system which is analogous to Theorem \ref{th2.6}.
\begin{t1} Let $X_{k:n}$ be the $k^{th}$ order statistic in a set of iid random variables from $F$. Then
\begin{itemize}
  \item $\overline\varepsilon J(X_{k:n};t)\geqslant\overline\varepsilon J(X_{k-1:n};t),$
  \item $\overline\varepsilon J(X_{k:n};t)\geqslant\overline\varepsilon J(X_{k:n+1};t),$
  \item $\overline\varepsilon J(X_{k:n};t)\geqslant\overline\varepsilon J(X_{k-1:n-1};t).$
\end{itemize}
\end{t1}
\hspace*{.2in} To end this section, we built characterization results based on DCPEx of $X_{n:n}$. The general characterization problem is to obtain whether $\overline\varepsilon J(X_{n:n};t)$ uniquely determines the underlying distribution function. We consider the following characterization result in analogy with Theorem \ref{th3.5}.
\begin{t1} Let $X$ and $Y$ be nonnegative absolutely continuous random variables with distribution functions $F(x)$ and $G(x)$, respectively. Then $F$ and $G$ belong to the same family of distributions, but for a change in location, if and only if
\begin{eqnarray*}\overline\varepsilon J(X_{n:n};t)=\overline\varepsilon J(Y_{n:n};t),\end{eqnarray*}
for all $n=n_j,~j\geqslant1$ such that $\sum_{j=1}^\infty n_j^{-1}=+\infty.$
\end{t1}
\hspace*{.2in} In the following theorem we characterize power distribution. This distribution includes the uniform distribution as a particular case.
\begin{t1} Let $X$ be a random variable with finite support $[0,b]$. Then $X$ has power distribution with $F(x)=(x/b)^c$ for $0<x<b$ and $c>0$, if and only if
\begin{eqnarray}\label{eq3.7}
\overline\varepsilon J(X_{n:n};t)=k m_F(t),
\end{eqnarray}
where $k$ is a constant.
\end{t1}
Proof: The `only if part' is straightforward with $k=\frac{-1}{2}\left(\frac{c+1}{2nc+1}\right)$. To prove the converse let us assume that (\ref{eq3.7}) holds. Then
$$\int_0^tF^{2n}(x)dx=2k(F(t))^{2n-1}\int_0^tF(x)dx,$$
which on differentiating with respect to $t$, gives
$$\phi_F(t)m_F(t)=\frac{1-2k}{2k(2n-1)}=k_1,~{\rm say}$$
Now, on using the relation $\phi_F(t)m_F(t)=1-m_F'(t)$, we obtain $m_F'(t)=1-k_1$ and using that $m_F(0)=0$ we have $m_F(t)=(1-k_1)t$ for $0<t<b$. Hence, using the inversion formula for $m_F(t)$, we obtain the stated distribution function.
\begin{c1} For an absolutely continuous random variable $X$ with finite support $[0,b]$
$$\overline\varepsilon J'(X_{n:n};t)=k,$$
where $k$ is a constant, if and only if $X$ follows power distribution as given in the above theorem.
\end{c1}

\section{Conclusions}
In recent years, there has been a great interest in the study of CREn and CPEn as an alternative measure of uncertainty. The basic idea is to replace the density function by survival/distribution function in Shannon's entropy. One important complementary dual measure of Shannon entropy is extropy which plays a vital role in the scoring of forecasting distributions, risk measure and independence. Here we consider cumulative analogue of extropy similar to CREn and CPEn and study it in context with extreme order statistics, i.e., $X_{1:n}$ and $X_{n:n}$. As the order statistic $X_{1:1}$ contains information about the location of the distribution on the real line. The results presented here generalize and enhance the related existing results in context with CREx and that are developed for order statistics. They are expected to be useful to the reliability theorists and the engineers, where entropy plays a vital role. Furthermore, this article is a first step in the study of CPEx also. Of course, other properties of CPEx are still waiting to be discovered.

%\section*{Acknowledgements}
%The financial support (Ref. No. SR/FTP/MS-016/2012) rendered by the Department of Science and Technology, Government of India, is acknowledged with thanks for carrying out this research work.

\end{document}